\documentstyle[12pt,righttag]{amsart}
\makeatletter

\newtheorem{nothm}{Theorem} 

\newtheorem{lem}{Lemma}

\theoremstyle{definition}

\theoremstyle{remark}

\newtheorem{ack}{Acknowledgment} 

\newcommand{\bbN}{\Bbb{N}}
\newcommand{\bbR}{\Bbb{R}}

\newcommand{\gam}{\gamma}
\newcommand{\Ome}{\Omega}
\newcommand{\al}{\alpha}
\newcommand{\del}{\delta}

\newcommand{\Gam}{\Gamma}

\newcommand{\supp}{\operatorname{supp}}

\newcommand{\lb}{\label}

\newcommand{\bs}{\backslash}
\newcommand{\bi}{\bibitem}
\newcommand{\emt}{\emptyset}

\newcommand{\lng}{\langle}
\newcommand{\rng}{\rangle}

\makeatletter
\def\@currentlabel{2.1}\label{e:dispaa}
\def\@currentlabel{2.21}\label{e:dispau}
\def\@currentlabel{2.22}\label{e:dispav}
\def\@currentlabel{2.23}\label{e:dispaw}
\def\@currentlabel{2.24}\label{e:dispax}
\def\theequation{\thesection.\@arabic\c@equation}

\makeatother

\makeatletter
\def\alphenumi{%
  \def\theenumi{\alph{enumi}}%
  \def\p@enumi{\theenumi}%
  \def\labelenumi{(\@alph\c@enumi)}}
\makeatother

\begin{document}
\title{Schoenberg's Problem on Positive Definite Functions}
\author{Alexander Koldobsky}
\maketitle
\bigskip\bigskip

\centerline{\sc University of Missouri--Columbia}
\centerline{\sc Department of Mathematics}
\centerline{\sc Columbia, MO 65211}

\bigskip\bigskip\bigskip
\section{Introduction}

In 1938, I.~J.~Schoenberg \cite{27} posed the following problem: For
which numbers
$\beta>0$ is the function $\exp (-\|x\|_q^\beta)$ positive definite on
$\bbR^n$?  Here $q>2$ and $\|x\|_q=(|x_1|^q+\cdots +|x_n|^q)^{1/q}$.
Denote by $B_n(q)$ the set of such numbers $\beta$.

We prove in this article that the functions $\exp (-\|x\|_q^\beta)$ are
not positive definite for all $n\geq 3$, $q>2$ and $\beta>0$, i.e.
$B_n(q)=\emt$ for every $q>2$ and $n\geq 3$.

Besides, $B_2(q)=(0,1]$ for every $q>2$.  In the case $n=2$ we have to
prove only that $B_2(q)\cap(1,\infty)=\emt$.  In fact, it is well--known
that $\exp(-\|x\|^\beta)$ is a positive definite function for every
two--dimensional norm and every $\beta\in (0,1]$, see Ferguson \cite{8},
Hertz \cite{11}, Lindenstrauss and Tzafriri \cite{19}, Dor \cite{6},
Misiewicz and Ryll--Nardziewski \cite{22}, Yost \cite{30}, Koldobsky
\cite{13} for different proofs.

Thus we give a complete answer to Schoenberg's question.

The case $q\in (0,2]$ was settled by Schoenberg \cite{27}.  Here one has
$B_n(q)=(0,q]$ for every $n\geq 2$.

Let us mention some facts which are valid for all $q>0$, see Schoenberg
\cite{27}:
\newcounter{aa}
\begin{list}{(\roman{aa})}{\usecounter{aa}}
\item $B_1(q)=(0,2]$;
\item $B_n(q)\supset  B_{n+1}(q)$ for every $n\in \bbN$;
\item if $\beta\in B_n(q)$ then $(0,\beta]\subset B_n(q)$.
\end{list}

Thanks to (ii) and (iii) it is enough to prove that $B_3(q)\cap
(0,2)=\emt$ to verify that $B_n(q)=\emt$ for every $n\geq 3$.

However, we treat here a more general situation.  Denote by $\phi_n(q)$
the class of all even functions $f:\bbR\to\bbR$ which are such that
$f(\|x\|_q)$ is a characteristic function, i.e. there exists a
probability measure $\mu$ on $\bbR^n$ with $\hat \mu(x)=f(\|x\|_q)$ for
all $x\in \bbR^n$.

It is clear that $\phi_n(q)\supset \phi_{n+1}(q)$.  If $f\in \phi_n(q)$
then by Bochner's theorem the function $f(\|x\|_q)$ is positive definite,
as well as $f$ itself.  Hence there exists a probability measure $\nu$ on
$\bbR$ with $\hat \nu=f$.

We prove the following

\begin{nothm}
Let $q>2$, $n\in\bbN$, $f\in\phi_n(q)$, $f\not\equiv 1$ and $\hat
\nu=f$.  Then in every one of the cases:
\newcounter{bb}
\begin{list}{\alph{bb})}{\usecounter{bb}}
\item $n\geq 3$, $\beta\in (0,2)$;
\item $n=2$, $\beta\in(1,2)$;
\item $n\geq 4$, $\beta\in (-1,0)$
\end{list}
the $\beta$--th moment of the measure $\nu$ is infinite, i.e.
$\int_{\bbR}|t|^\beta d\nu(t)=\infty$.
\end{nothm}

Using a) one can immediately get an answer to Schoenberg's question for
$n\geq 3$.  Indeed, for every $\beta\in (0,2)$, $\exp(-|t|^\beta)$ is the
characteristic function of the $\beta$--stable measure on $\bbR$.  This
measure has finite moments of all positive orders less than $\beta$,  see
Zolotarev \cite{29} or the text preceding Lemma~\ref{l4} below.  By a)
the function $\exp(-|t|^\beta)\not\in \phi_3(q)$ for every $q>2$.  By
Bochner's theorem the function $\exp(-\|x\|_q^\beta)$ is not positive
definite on $\bbR^3$.  Thus $B_3(q)\cap (0,2)=\emt$ for every $q>2$ and
we are done.

For the same reason, it follows from b) that $B_2(q)\cap (1,2)=\emt$ for
every $q>2$.

The statements a) and c) show that the measure $\nu$ corresponding to a
function $f\in \phi_n(q)$, $n\geq 4$, $q>2$ must have a very special
behavior at infinity and at zero.  Nevertheless, the following question
is open: For $n\geq 3$ and $q>2$, are there any functions in the class
$\phi_n(q)$ besides the function $f\equiv 1$?

The classes $\phi_n(q)$ have been investigated by several authors.
Schoenberg \cite{28} described the classes $\phi_n(2)$ and
$\phi_\infty(2)=\bigcap_n\phi_n(2)$ completely.  A similar result was
obtained for the classes $\phi_n(1)$ by Cambanis et al \cite{4}.
Bretagnolle et al \cite{3} described the classes $\phi_\infty(q)$ for all
$q>0$.  In particular, for very $q>2$ the class $\phi_\infty(q)$ contains
no functions besides $f\equiv 1$.  For some partial results on the
classes $\phi_n(q)$, $0<q<2$, see Richards \cite{25,26} and Misiewicz
\cite{20}.  Misiewicz \cite{21} proved that for $n\geq 3$ a function
$f(\max(|x_1|, \ldots, |x_n|))$ is positive definite only if $f\equiv 1$.
This result gives an answer to Problem 2 from Schoenberg \cite{27}.
References related to the topic include also Askey \cite{1}, Berg and
Ressel \cite{2}, Kuelbs \cite{15}, Eaton \cite{7}, Christensen and Ressel
\cite{5}, Kuritsyn \cite{16}, Kuritsyn and Shestakov \cite{17}, Misiewicz
and Scheffer \cite{23}, Aharoni et al \cite{31}.

The well known theorem of Bretagnolle et al \cite{3} states that a Banach
space $(E, \|\cdot\|)$ is isometric to a subspace of $L_\beta([0,1])$,
$\beta\in[1,2]$ iff the function $\exp (-\|x\|^\beta)$ is positive
definite.  Dor \cite{6} pointed out all the pairs of numbers $\beta$,
$q\in[1,\infty)$ for which the space $l_q^n$, $n\geq 2$, is isometric to
a subspace of $L_\beta([0,1])$.  Combining these results one can prove
that $B_n(q)\cap (1,\infty)=\emt$ for every $n\geq 2$ and $q>2$.  So
Schoenberg's problem was open only for $\beta\in (0,1)$.

In spite of the connection between Schoenberg's problem and isometries,
in this paper we don't deal with isometries directly.  Our main tool is
the Fourier transform of distributions.

\bigskip
\section{Some Applications of the Fourier Transform}

Let $(E, \|\cdot\|)$ be an $n$--dimensional Banach space and
$f:\bbR\to\bbR$ be an even function such that there exists a probability
measure $\mu$ on $\bbR^n$ with $\hat\mu(x)=f(\|x\|)$ for all
$x\in\bbR^n$.  Measures with characteristic functions of the form
$f(\|x\|)$ are called $E$--stable or pseudo--isotropic, since all
one--dimensional projectional of such measures are equal up to a scale
parameter.  We give an easy proof of this fact.

\begin{lem}[see Levy \cite{18}, Eaton \cite{7}]\lb{l1}
If $\mu$ is a probability measure on $\bbR^n$ with the characteristic
function $f(\|x\|)$ then:
\begin{list}{\alph{bb})}{\usecounter{bb}}
\item there exists a probability measure $\nu$ on $\bbR$ with $\hat
\nu=f$;
\item for every $x\in\bbR^n$, $x\not= 0$, the image of the measure $\mu$
under the mapping $\xi\to\lng x,\xi \rng / \|x\|$ from $\bbR^n$ to $\bbR$
coincides with the measure $\nu$ (here $\lng x,\xi\rng$ stands for the
scalar product).
\end{list}
\end{lem}

\begin{pf}
Fix an element $x\in \bbR^n$, $x\not= 0$.  Let $\nu$ be the image of the
measure $\mu$ under the mapping $\xi\to \lng x,\xi\rng / \|x\|$.  Put
$y=\lng x,\xi\rng /\|x\|$.  Then for every $k\in \bbR$
\begin{align*}
& f(k\|x\|)=\hat\mu(kx)=\int_{\bbR^n} \exp (-i\lng kx, \xi\rng )\,
d\mu
(\xi)\\
&\quad = \int_{\bbR^n}\exp (-ik\|x\|(\lng x, \xi\rng / \|x\|))\, d\mu
(\xi)=\int_{\bbR} \exp (-ik\|x\|y)\, d\nu (y)=\hat \nu(k\|x\|).
\end{align*}
Hence $\hat \nu =f$ and $\nu$ doesn't depend on the choice of $x\in
\bbR^n$, $x\not= 0$.
\end{pf}

Let $S(\bbR^n)$ be, as usual, the space of rapidly decreasing infinitely
differentiable functions.  Denote by $S'(\bbR^n)$ the space of
distribution over $S(\bbR^n)$.

Let $\Ome$ be an open subset of $\bbR^n$.  A distribution $g\in
S'(\bbR^n)$ is called positive (negative) on $\Ome$ if $\lng g, \psi\rng
\geq 0$ $(\lng g, \psi\rng \leq 0)$ for every non--negative function
$\psi\in S(\bbR^n)$ with $\supp \psi \subset \Ome$.

\begin{lem}[cf. Koldobsky \cite{12}]\lb{l2}
Let $\beta\in (-1,\infty)$, $\beta\not= 0,2,4,\ldots$.  Let $\psi\in
S(\bbR^n)$ be a function with $0\not\in \supp \psi$.  Then for every
$\xi\in\bbR^n$, $\xi\not= 0$,
$$
\int_{\bbR^n} |\lng x,\xi \rng|^\beta \hat \psi (x)\, dx=c_\beta
\int_{\bbR} |t|^{-1-\beta} \psi(t\xi)\, dt,\;
c_\beta=\frac{2^{\beta+1}\pi^{1/2}\Gam ((\beta+1)/2)}{\Gam(-\beta/2)}.
$$
\end{lem}

\begin{pf}
It is well--known that $(|x|^\beta)^\wedge (t)=c_\beta|t|^{-1-\beta}$,
$t\not= 0$, for all $\beta\in (-1,\infty)$, $\beta\not= 0,2,4,\ldots$,
see Gelfand and Shilov \cite{9}.  By the Fubini theorem
\begin{equation}\lb{2.1}
\int_{\bbR^n}|\lng x,\xi\rng|^\beta \hat\psi(x)\, dx=\int_{\bbR}|z|^\beta
\Big(\int_{\lng x,\xi\rng=z} \hat \psi (x)\, dx\Big)\, dz=\Big\lng
|z|^\beta,
\int_{\lng x, \xi \rng =z} \hat \psi (x)\, dx\Big\rng.
\end{equation}
The function $t\to 2\pi \psi(-t\xi)$ is the Fourier transform of the
function $z\to \int_{\lng x, \xi \rng=z} \hat\psi (x)\, dx$ (it is a
simple property of the Radon transform, see Gelfand et al \cite{10}).
Therefore, we can continue the equality (\ref{2.1}):
$$
=c_\beta \lng |t|^{-1-\beta}, \psi(t\xi)\rng=c_\beta \int_{\bbR}
|t|^{-1-\beta}\psi(t\xi)\, dt.\hspace*{2mm} \qed
$$
\renewcommand{\qed}{}\end{pf}

\begin{lem}\lb{l3}
Let $f$, $\mu$ and $\nu$ be as in Lemma~\ref{l1}.  Additionally, assume
that $f\not\equiv 1$ and $\int_{\bbR}|t|^\beta \, d\nu (t)<\infty$ for a
number $\beta\in (-1,2)$, $\beta\not= 0$.  Then, if $\beta\in (-1,0)$ the
distribution $(\|x\|^\beta)^\wedge$ is positive on $\bbR^n \bs \{0\}$.
If $\beta\in (0,2)$ the distribution $(\|x\|^\beta)^\wedge$ is negative
on $\bbR^n \bs \{0\}$.
\end{lem}

\begin{pf}
Since $f \not\equiv 1$, the measure $\nu$ is not supported at zero and
$\int_{\bbR}|t|^\beta \, d\nu (t)>0$.  By Lemma~\ref{l1} we have
$$
\int_{\bbR^n}|\lng x,\xi\rng|^\beta \, d\mu (\xi)=\|x\|^\beta
\int_{\bbR^n}\left|\frac{\lng x,\xi\rng}{\|x\|}\right|^\beta\, d\mu
(\xi)=\|x\|^\beta \int_{\bbR}|t|^\beta\, d\nu (t)
$$
for every $x\in \bbR^n$.  It follows from Lemma~\ref{l2} and the Fubini
theorem that
\begin{align*}
&\lng(\|x\|^\beta)^\wedge, \psi\rng=\lng\|x\|^\beta,
\hat\psi\rng=\dfrac1{\int_{\bbR} |t|^\beta \, d\nu (t)} \int_{\bbR^n}\,
d\mu (\xi)\int_{\bbR^n} |\lng x,\xi\rng|^\beta \hat\psi (x)\, dx\\
&\quad =\frac{c_\beta}{\int_{\bbR}|t|^\beta \, d\nu (t)}\int_{\bbR} \,
d\mu
(\xi) \int_{\bbR} |t|^{-1-\beta}\psi(t\xi)\, dt
\end{align*}
for every function $\psi\in S(\bbR^n)$ with $0\not\in \supp \psi$.  If
the function $\psi$ is non--negative, the quantity in the right--hand
side has the same sign as the number $c_\beta$.  It suffices to remark
that $c_\beta>0$ if $\beta\in (-1,0)$ and $c_\beta<0$ if $\beta\in
(0,2)$.
\end{pf}

For arbitrary $q>0$, define the function $\gam_q$ on $\bbR$ by
$\gam_q(t)=(\exp (-|x|^q))^\wedge (t)$, $t\in \bbR$.  Then
$$
\lim\limits_{t\to\infty} t^{1+q}\gam_q(t)=2\Gam (q+1) \sin (\pi q/2),
$$
see Polya and Szego \cite{24}, Part 3, Problem 154.  Therefore, the
integral
$$
S_q(\al)=\int_{\bbR}|t|^\al \gam_q(t)\, dt
$$
converges absolutely for every $\al\in (-1,q)$.

\begin{lem}[for the case $q\in (0,2)$, see Zolotarev \cite{29} or
Koldobsky \cite{12}]
Let\\ $q>2$.  Then for every $\al\in (-1,q)$, $\al\not=
0,2,\ldots,
2[q/2]$,
$$
S_q(\al)=2^{\al+2}\pi^{1/2} \Gam(-\al / q)\Gam ((\al+1)/2)/ (q\Gam
(-\al/2)).
$$
In particular, $S_q(\al)>0$ if $\al\in (-1,2)$, $\al\not= 0$, and
$S_q(\al)<0$ if $\al\in (2, \min (4,q))$.
\lb{l4}
\end{lem}

\begin{pf}
Assume $-1<\al<0$.  By the Parseval theorem
\begin{align}
\begin{split}
& S_q(\al)=\int_{\bbR} |t|^\al
\gam_q(t)\, dt=c_\al \int_{\bbR} |z|^{-1-\al}\exp (-|z|^q)\, dz\\
&\quad =\frac{2^{\al+1}\pi^{1/2}\Gam ((\al+1)/2)}{\Gam(-\al/2)}\cdot
\frac{2\Gam(-\al/q)}{q}.
\lb{2.2}
\end{split}
\end{align}
If we allow $\al$ to assume complex values then all functions of $\al$ in
(\ref{2.2}) are analytic in the domain $\{-1< R e \al <q, \al\not=
0,2,\ldots,2[q/2]\}$.  Since analytic continuation from the interval
$(-1,0)$ is unique, the equality (\ref{2.2}) remains valid for all
$\al\in(-1,q)$, $\al\not= 0,2,\ldots,2[q/2]$.  To complete the proof,
note that $\Gam (z)>0$ if $z>0$ or $z\in (-2,-1)$, and $\Gam(z)<0$ if
$z\in (-1,0)$.
\end{pf}

Let us compute the Fourier transform of the function $\|x\|_q^\beta$.

\begin{lem}\lb{l5}
Let $q>0$, $n\in\bbN$, $-n<\beta<qn$, $\beta/q \not\in \bbN\cup \{0\}$,
$\xi=(\xi_1, \ldots, \xi_n)\in \bbR^n$, $\xi_k\not= 0$, $1\leq k\leq n$.
Then
\begin{equation}\lb{2.3}
((|x_1|^q+\cdots +|x_n|^q)^{\beta/q})^\wedge
(\xi)=\frac{q}{\Gam(-\beta/q)} \int_0^\infty t^{n+\beta-1}\prod_{k=1}^n
\gam_q (t\xi_k)\, dt
\end{equation}
\end{lem}

\begin{pf}
Assume $-1<\beta<0$.  By definition of the $\Gam$--function
$$
(|x_1|^q+\cdots +|x_n|^q)^{\beta /q}
=\frac{q}{\Gam(-\beta/q)}\int_0^\infty y^{-1-\beta} \exp
(-y^q(|x_1|^q+\cdots + |x_n|^q))\, dy.
$$
For every fixed $y>0$,
$$
(\exp (-y^q(|x_1|^q+\cdots +|x_n|^q))^\wedge (\xi) =y^{-n}\prod_{k=1}^n
\gam_q (t\xi_k).
$$
Put $t=1/y$.  We have
\begin{align*}
&((|x_1|^q+\cdots+|x_n|^q)^{\beta/q})^\wedge
(\xi)=\frac{q}{\Gam(-\beta/q)}\int_0^\infty y^{-n-\beta-1} \prod_{k=1}^n
\gam_q (\xi_k /y)\, dy\\
&\quad = \frac{q}{\Gam(-\beta/ q)}\int_0^\infty t^{t+\beta-1}
\prod_{k=1}^n \gam_q (t\xi_k)\, dt.
\end{align*}

The latter integral converges if $-n<\beta<qn$ because the function
$t\to \prod_{k=1}^n \gam_q (t\xi_k)$ decreases at infinity like
$t^{-n-nq}$ (remind that $\xi_k\not= 0$, $1\leq k\leq n$; see the text
preceding Lemma~\ref{l4}).

If $\beta$ is allowed to assume complex values then the both sides of
(\ref{2.3}) are analytic functions of $\beta$ in the domain $\{-n<Re
\beta <nq, \beta/q \not\in \bbN\cup \{0\}\}$.  These two functions admit
unique analytic continuation from the interval $(-1,0)$.  Thus the
equality (\ref{2.3}) remains valid for all $\beta\in (-n,qn)$,
$\beta/q\not\in \bbN\cup \{0\}$ (see Gelfand and Shilov \cite{9} for
details of analytic continuation in such situations).
\end{pf}

For $q\in [1,2]$, Lemma~\ref{l5} was proved in Koldobsky \cite{14}.

\bigskip
\section{Proof of Theorem}

Let $n\in\bbN$, $n\geq 2$, $-n<\beta<qn$, $q>2$.  Consider the integral
\setcounter{equation}{0}
\begin{align}
\begin{split}
& J_n(\al_1,\ldots,\al_{n-1})\\
& \quad =\int_{\bbR^{n-1}} |\xi_1|^{\al_1}\cdot \cdots \cdot
|\xi_{n-1}|^{\al_{n-1}}\cdot \Big(\int_0^\infty t^{n+\beta-1}
\gam_q(t)\prod_{k=1}^{n-1} \gam_q(t\xi_k)\, dt\Big)\,
d\xi_1\cdot\cdots\cdot d\xi_{n-1}\\
&\quad = \int_0^\infty t^{n+\beta-1} \gam_q(t)\Big(\prod_{k=1}^{n-1}
\int_{\bbR}|\xi_k|^{\al_k}\gam_q (t\xi_k)\, d\xi_k\Big)\, dt\\
&\quad = \int_0^\infty t^{-\al_1-\cdots -\al_{n-1}+\beta}\gam_q(t)\,
dt\cdot
\prod_{k=1}^{n-1} \int_{\bbR} |z|^{\al_k}\gam_q (z)\, dz\\
&\quad =S_q (\al_1)\cdot S_q(\al_2)\cdot \cdots\cdot S_q
(-\al_1-\al_2-\cdots
-\al_{n-1}+\beta).
\lb{3.1}
\end{split}
\end{align}

If all the numbers $\al_1,\ldots, \al_{n-1}$,
$-\al_1-\cdots-\al_{n-1}+\beta$ belong to the interval $(-1,q)$ then the
integrals in (\ref{3.1}) converge absolutely.  Therefore, the Fibini
theorem is applicable and the integral $J_n (\al_1,\ldots,\al_{n-1})$
converges.

Case a): $n\geq 3$, $\beta\in(0,2)$.

Note that $\phi_n(q)\supset \phi_{n+1}(q)$, so it is enough to prove the
theorem for $n=3$.

Let $n=3$ and assume that there exists a function $f\in \phi_3(q)$,
$f\not\equiv 1$, with $\int_{\bbR}|t|^\beta \, d\nu (t)<\infty$, where
$\hat\nu =f$ and $\beta\in (0,2)$.

Put $\al_1=\al_2=-1+\del$ with $\del\in (0,1)$ and
$(2+\beta-\min(4,q))<\del<\beta/2$.  Then $\al_1, \al_2\in(-1,0)$ and
$-\al_1-\al_2+\beta\in(2,\min(4,q))$.  By Lemma~\ref{l4},
$J_3(\al_1,\al_2)=S_q(\al_1)\cdot S_q(\al_2)\cdot
S_q(-\al_1-\al_2+\beta)<0$.

On the other hand, by Lemma~\ref{l3} the function
$(\|x\|_q^\beta)^\wedge$ is negative on $\bbR^n \bs \{0\}$.  Since
$\Gam(-\beta/ q)<0$, it follows from Lemma~\ref{l5} that
$$
h(\xi_1, \xi_2)=\int_0^\infty t^{\beta+2}
\gam_q(t\xi_1)\gam_q(t\xi_2)\gam_q(t)\, dt\geq 0
$$
for every $\xi_1, \xi_2 \not= 0$.  Hence
$$
J_3(\al_1,\al_2)=\int_{\bbR^2}|\xi_1|^{\al_1}\cdot |\xi_2|^{\al_2}
h(\xi_1,\xi_2)\, d\xi_1 \, d\xi_2>0
$$
and we get a contradiction.

In the cases b) and c) proofs are similar.  In the case b) we put
$\al_1=-1+\del$ with $\del\in (0,1)$ and
$1+\beta-\min(4,q)<\del<\beta-1$.  In the case c) we may restrict
ourselves to $n=4$ and put $\al_1=\al_2=\al_3=-1+\del$ with
$\del\in(0,1)$ and $(3+\beta-\min(4,q))/3<\del <(\beta+1)/3$.  Then the
numbers $J_2(\al_1)$ and $J_4(\al_1, \al_2, \al_3)$ are negative.  On the
other hand, if there exists a function $f$ for which the measure $\nu$
has finite $\beta$--th moment, then Lemma~\ref{l3} and Lemma~\ref{l5}
imply $J_2(\al_1)\geq 0$ and $J_4(\al_1, \al_2, \al_3)\geq 0$ and we get
a contradiction (note that $\Gam(-\beta/q)>0$ if $\beta\in (-1,0)$ and
$\Gam(-\beta/q)<0$ if $\beta\in(1,2)$.)

\begin{ack}
I am grateful to Jolanta Misiewicz for helpful discussions during my stay
in Technical University of Wroclaw in 1989.
\end{ack}

\end{document}